\documentclass[11pt,a4paper]{amsart}
\usepackage[utf8]{inputenc}
\usepackage[T1]{fontenc}
\usepackage{amsmath,amssymb,amsthm,mathtools}
\usepackage{enumitem}
\usepackage[margin=2.8cm]{geometry}
\usepackage{xcolor}
\usepackage[colorlinks=true,linkcolor=blue!60!black,citecolor=green!40!black,urlcolor=blue!60!black]{hyperref}

\newif\ifnotes
\notestrue

\definecolor{revblue}{rgb}{0,0.25,0.8}
\definecolor{revred}{rgb}{0.78,0,0}
\definecolor{revgreen}{rgb}{0,0.45,0.15}
\definecolor{revviolet}{rgb}{0.5,0,0.6}
\newcommand{\corr}[1]{\ifnotes{\color{revblue}\relax#1}\else#1\fi}
\newcommand{\sug}[2]{\ifnotes{\color{revred}\,[\textsf{\small Sugerencia:}~#1\ifx\relax#2\relax\else\ {\sffamily\small(#2)}\fi]}\fi}

\newtheorem{maintheorem}{Theorem}

\newtheorem{theorem}{Theorem}[section]
\newtheorem{lemma}[theorem]{Lemma}
\newtheorem{proposition}[theorem]{Proposition}
\newtheorem{corollary}[theorem]{Corollary}
\theoremstyle{definition}
\newtheorem{definition}[theorem]{Definition}
\newtheorem{example}[theorem]{Example}
\newtheorem{question}[theorem]{Question}
\theoremstyle{remark}
\newtheorem{remark}[theorem]{Remark}

\newcommand{\N}{\mathbb{N}}
\newcommand{\Lang}{\mathcal{L}}
\newcommand{\G}{\mathcal{G}}
\newcommand{\Cp}{\mathcal{C}^{p}}
\newcommand{\Cs}{\mathcal{C}^{s}}
\newcommand{\eps}{\varepsilon}
\newcommand{\nz}[1]{|#1|_{0}}
\newcommand{\no}[1]{|#1|_{1}}
\newcommand{\Fol}{F}

\begin{document}

\title[Frequency shifts]{Frequency shifts: non-soficity, absence of specification\\ and intrinsic ergodicity}

\author{Carlos Reyes}
\address{Instituto de Matem\'atica, Estat\'istica e Ci\^encia da Computa\c{c}\~ao, Universidade de S\~ao Paulo, Brazil}
\email{carlosdhers123@gmail.com}

\date{\today}

\begin{abstract}
For $p_c\in(0,1)$, the frequency shift $\Sigma(p_c)\subset\{0,1\}^{\N}$, introduced by Becker and Baraviera, is the closure of the shift orbits of the set of sequences in which a symbol $1$ is followed by $1$ exactly when the density of ones up to that position does not exceed $p_c$. We show that $\Sigma(p_c)$ is never sofic and never has the (weak) specification property. For $p_c=\frac{k}{k+1}$, with $k\ge1$, we describe $\Sigma(p_c)$ completely: $\Sigma(\frac12)$ is the coded system generated by $\{0^m1^m\}_{m\ge1}$, while for $k\ge2$ the shift $\Sigma(\frac{k}{k+1})$ is not transitive, but every invariant probability measure is supported on the coded system generated by $\{0^m1^{km}\}_{m\ge1}$. For every $k \ge 1$, the topological entropy of $\Sigma(\frac{k}{k+1})$ is $\frac{\log2}{k+1}$, and, using the Climenhaga--Thompson criterion, we prove that $\Sigma(\frac{k}{k+1})$ is intrinsically ergodic. 

\end{abstract}


\maketitle

%==================================================================
\section{Introduction}\label{sec:intro}
%==================================================================
Entropy is the natural quantity for describing how a system behaves at equilibrium. Given a dynamical system, an invariant measure whose entropy equals the topological entropy is called a \emph{measure of maximal entropy} (MME). Such measures are precisely the equilibrium states of the zero potential, and the question of their uniqueness is part of the thermodynamic formalism of Sinai, Ruelle and Bowen \cite{sinai1972gibbs,ruelle1978thermodynamic,bowen2008equilibrium}. When the MME is unique, the system is called \emph{intrinsically ergodic}; in the language of statistical mechanics, the coexistence of several equilibrium states is interpreted as a phase transition (see, e.g., \cite{haydn2013phase} for one-dimensional subshifts).

Throughout this paper\corr{,} we work with \emph{one-sided} subshifts, that is, closed subsets $X\subseteq\{0,1\}^{\N}$ with $\sigma(X)\subseteq X$, where $\sigma$ is the left shift, which is not invertible. This is the natural setting for the systems studied here, since their defining rule involves the density of ones counted from the first coordinate. Much of the classical theory is formulated for two-sided subshifts. Notions defined through the language, such as topological entropy, transitivity and specification, carry over verbatim, and the invariant measures of a one-sided subshift correspond bijectively, with the same entropies, to those of its natural extension, which is a two-sided subshift \cite[Prop.~2.1]{climenhaga2012intrinsic}. However, the language of a one-sided subshift may contain words that cannot be extended to the left, and such words can change its topological properties. As we shall see, this happens for some of the frequency shifts studied here, which fail to be transitive although all their invariant measures are supported on a transitive subsystem (Theorem~\ref{thm:C}(b)). When a classical result is stated only for two-sided subshifts or for homeomorphisms, we recall its one-sided version, as for the existence of measures of maximal entropy (Section~\ref{sec:prelim}), use a definition adapted to one-sided subshifts, as for sofic subshifts (Section~\ref{sec:nonsofic}), or recall that it remains valid in this setting, as for the Climenhaga--Thompson criterion (Section~\ref{sec:IE}).

For one-sided subshifts, as for two-sided ones, measures of maximal entropy always exist, because the entropy map $\mu\mapsto h_\mu(\sigma)$ is upper semicontinuous when the alphabet is finite, for a discussion about this point see \cite{Iommi2020}. The question is therefore one of uniqueness, and it can be attacked through the combinatorics of the language: is the language rich enough to force a unique MME, or does it allow several? The classical sufficient condition is Bowen's specification property \cite{bowen1971periodic,bowen1974some}, which, in the symbolic setting, says that any two words of the language can be joined by a connecting word of a fixed length.
It holds, for example, for every mixing subshift of finite type. Many natural subshifts do not satisfy it, and for those the problem is more delicate. Typical examples are the $\beta$-shifts, which have specification only for a Lebesgue-null set of parameters \cite{schmeling1997symbolic} but are nevertheless intrinsically ergodic \cite{hofbauer1978beta}, and the bounded density shifts introduced by Stanley \cite{stanley}.
Both of these families are defined by upper bounds, and both are therefore \emph{hereditary} \cite{kwietniak2011topological}: if $x\in X$ and $y\le x$ coordinatewise, then $y\in X$. The $\beta$-shifts are bounded from above by a sequence coming from a $\beta$-expansion, and the bounded density shifts are restricted by length-dependent bounds on the sums of letters in subwords. Stanley \cite{stanley} characterized when a bounded density shift is of finite type, is sofic or has the specification property, in terms that are strikingly similar to those known for $\beta$-shifts \cite{schmeling1997symbolic}.

Among the tools available for proving intrinsic ergodicity in the absence of specification, one of the most flexible is the criterion of Climenhaga and Thompson \cite{climenhaga2012intrinsic}, which only requires specification on part of the language. Roughly speaking, one writes each word of the language as a concatenation $uvw$ and asks that
(I) the central words $v$ lie in a collection $\mathcal G$ with specification;
(II) the collections of possible prefixes $u$ and suffixes $w$ have entropy strictly smaller than that of the whole system; and
(III) whenever $|u|,|w|\le M$, the word $uvw$ can be extended to an element of $\mathcal G$ by attaching at most $\tau(M)$ symbols on each side.
Under these hypotheses, the system has a unique measure of maximal entropy (Theorem~\ref{thm:CT}). For $\beta$-shifts, the criterion gives intrinsic ergodicity with a very short argument, and it was also the main tool in \cite{garcia-ramos_pavlov_reyes2024}, where the authors obtained a fairly general sufficient condition for the intrinsic ergodicity of bounded density shifts. 

In this paper we apply this strategy to a family of one-sided subshifts on $\{0,1\}$ introduced by Becker and Baraviera \cite{becker2023symbolic}. Fix $p_c\in[0,1]$. The set $P(p_c)\subset\{0,1\}^{\N}$ consists of the sequences $x$ such that, whenever $x_i=1$, the next symbol is
\[
x_{i+1}=1 \ \text{ if } \ \delta_1(x_1\cdots x_i)\le p_c,
\qquad
x_{i+1}=0 \ \text{ if } \ \delta_1(x_1\cdots x_i)> p_c,
\]
where $\delta_1$ denotes the density of ones in a word; after a $0$, both symbols are allowed. The \emph{frequency shift} is $\Sigma(p_c)=\overline{\bigcup_{n\ge0}\sigma^n(P(p_c))}$. In \cite{becker2023symbolic}, it is shown that $\Sigma(0)$ and $\Sigma(1)$ are subshifts of finite type, that $\Sigma(p_c)$ is not of finite type for $p_c\in(0,1)$, and that $\Sigma(p_c)$ has positive topological entropy for $p_c\in[0,1)$.

These frequency shifts are also defined through densities, but in a different way: the density of ones in a prefix does not bound the sequence from above; it \emph{decides} the next symbol after each $1$. As a consequence, for $p_c\in(0,1)$, they are not hereditary (Proposition~\ref{prop:not-hereditary}); moreover, none of them is sofic and none has the specification property (Theorems~\ref{thm:A} and~\ref{thm:B}). Still, for the parameters $p_c=\frac{k}{k+1}$, $k\ge1$, we prove that they are intrinsically ergodic, by identifying inside them a coded system that carries every invariant measure.

For this coded system, we use a variant of the natural decomposition of its language \cite{climenhaga2012intrinsic}, in which words are cut into a suffix of a generator, a concatenation of generators and a prefix of a generator. The natural decomposition does not cover the words lying strictly inside a single generator; these form a family of zero entropy, and we add them to the collection of suffixes (Remark~\ref{rem:natural}). This enlarged decomposition is adapted to our examples, and we do not claim that it works for coded systems in general.

The aim of this paper is to locate these systems within the usual hierarchy of symbolic systems---sofic shifts, shifts with specification, coded and synchronized systems---and to address intrinsic ergodicity. Our main tool is an elementary reformulation of the transition rule. Put $\beta=p_c/(1-p_c)$ (this $\beta$ is unrelated to the $\beta$-shifts mentioned above; we keep the notation of \cite{becker2023symbolic}) and, for a word $w$, define the \emph{discrepancy}
\[
D(w)=\beta\,\nz{w}-\no{w},
\]
where $\nz{w}$ and $\no{w}$ denote the numbers of zeros and ones in $w$. Then a $1$ is followed by $1$ if and only if the discrepancy of the prefix ending at it is non-negative (Lemma~\ref{lem:rule}). Each $0$ raises the discrepancy by $\beta$ and each $1$ lowers it by $1$, and the length of each run of ones is determined by the discrepancy at the start of the run (Lemma~\ref{lem:runs}). The main results are the following:
\begin{maintheorem}\label{thm:A}
For every $p_c\in(0,1)$, the subshift $\Sigma(p_c)$ is not sofic.
\end{maintheorem}

\begin{maintheorem}\label{thm:B}
For every $p_c\in(0,1)$, the subshift $\Sigma(p_c)$ does not have the weak specification property.
\end{maintheorem}

For $k\ge1$, let $G_k=\{0^m1^{km}:m\ge1\}$ and let $X_k$ be the coded system generated by $G_k$.

\begin{maintheorem}\label{thm:C}
Let $k\ge1$.
\begin{enumerate}[label=\textup{(\alph*)}]
\item If $k=1$, then $\Sigma(\frac12)=X_1$; in particular, $\Sigma(\frac12)$ is coded, transitive and synchronized.
\item If $k\ge2$, then $X_k\subsetneq\Sigma(\frac{k}{k+1})$, the subshift $\Sigma(\frac{k}{k+1})$ is not transitive, and every $\sigma$-invariant Borel probability measure on $\Sigma(\frac{k}{k+1})$ is supported on $X_k$.
\item $h_{\mathrm{top}}(\Sigma(\frac{k}{k+1}))=\frac{\log2}{k+1}$.
\end{enumerate}
\end{maintheorem}

\begin{maintheorem}\label{thm:D}
For every $k\ge1$, the subshift $\Sigma(\frac{k}{k+1})$ has a unique measure of maximal entropy.
\end{maintheorem}

Theorem~\ref{thm:A} gives, in particular, a self-contained proof that, for every $p_c \in (0,1)$, $\Sigma(p_c)$ is not of finite type (Corollary~\ref{cor:notSFT}; compare Remark~\ref{rem:BB-thm1}). Theorems~\ref{thm:B} and~\ref{thm:D} show that the frequency shifts with $p_c=\frac{k}{k+1}$ fit into the framework of \emph{intrinsic ergodicity beyond specification} developed by Climenhaga and Thompson~\cite{climenhaga2012intrinsic}. 

\emph{Relation with other coded systems.} The coded systems $X_k$ are reminiscent of the $S$-gap shifts of \cite{climenhaga2012intrinsic}, whose generators $0^n1$ ($n\in S$) only constrain the gaps between consecutive ones; here, instead, the length of each block of ones is tied to the length of the preceding block of zeros. Moreover, $\Sigma(\frac12)=X_1$ is a natural example of a synchronized, non-sofic coded system arising from a density constraint.

\emph{Organization of the paper.} In Section~\ref{sec:prelim} we recall the basic notions and introduce the discrepancy, which governs the lengths of the runs of ones (Lemma~\ref{lem:runs}) and leads to the key estimate on anchored runs (Lemma~\ref{lem:anchored}). Theorems~\ref{thm:A} and~\ref{thm:B} are proved in Sections~\ref{sec:nonsofic} and~\ref{sec:nospec}, respectively. In Section~\ref{sec:k} we describe $\Sigma(\frac{k}{k+1})$ and prove parts~(a) and~(b) of Theorem~\ref{thm:C}; part~(c) and Theorem~\ref{thm:D} are proved in Section~\ref{sec:IE}, using the Climenhaga--Thompson criterion. Section~\ref{sec:questions} collects some open questions.
%==================================================================
\section{Preliminaries}\label{sec:prelim}
%==================================================================

\subsection{Words and subshifts}

A \emph{word} is a finite sequence $w=w_1\cdots w_n$ over $\{0,1\}$. We write $|w|=n$ for its length, $\eps$ for the empty word and $\nz{w},\no{w}$ for the numbers of zeros and ones in $w$. Concatenation is denoted by juxtaposition, and for a sequence $x=x_1x_2\cdots\in\{0,1\}^{\N}$,  we write $x_{[i,j]}=x_i\cdots x_j$, and similarly $w_{[i,j]}$ for a word $w$. For a word $u$ and $n\ge0$, $u^n$ denotes the concatenation of $n$ copies of $u$ (with $u^0=\eps$), and $u^\infty=uuu\cdots$.  A word $u$ is a \emph{factor} of $w$ (or of $x$) if $w=aub$ (respectively $x=aub\cdots$) for some words $a,b$. A \emph{maximal run of ones} in $w$ is a factor $w_{[i,j]}=1^{j-i+1}$ such that $i=1$ or $w_{i-1}=0$, and $j=|w|$ or $w_{j+1}=0$. We call it \emph{closed} if $j<|w|$. Maximal runs of zeros, and maximal runs in a sequence $x\in\{0,1\}^{\N}$, are defined in the same way (for a sequence, the alternative $j=|w|$ does not occur, so every maximal run that ends is closed).

Let $\sigma$ be the left shift map on $\{0,1\}^{\N}$. A \emph{subshift} is a closed set $X$ with $\sigma(X)\subseteq X$. Its \emph{language} $\Lang(X)$ is the set of all factors of points of $X$, and $\Lang_n(X)$ is the set of words of length $n$ in $\Lang(X)$. A subshift is determined by its language: $x\in X$ if and only if every factor of $x$ belongs to $\Lang(X)$. For $w\in\Lang(X)$, the \emph{follower set} of $w$ is $\Fol_X(w)=\{v:wv\in\Lang(X)\}$.

\begin{definition}\label{def:basic}
Let $X\subseteq\{0,1\}^{\N}$ be a subshift.

\begin{enumerate}[label=(\roman*)]
\item $X$ is a \emph{subshift of finite type} (SFT) if there is a finite set $\mathcal F$ of words, called \emph{forbidden words}, such that
\begin{equation*}
X=\{x\in\{0,1\}^{\mathbb N}:\ \text{no factor of } x \text{ belongs to } \mathcal F\}.
\end{equation*}

\item $X$ is \emph{transitive} if for all $u,v\in\Lang(X)$ there is a word $w$ with $uwv\in\Lang(X)$.
\item $X$ has the \emph{specification property} if there is $N\ge0$ such that for all $u,v\in\Lang(X)$ there is a word $w$ with $|w|=N$ and $uwv\in\Lang(X)$. For subshifts, this is the usual form of Bowen's specification property \cite{bowen1971periodic,bowen1974some}; see \cite{kwietniak2015panoramaspecificationlikepropertiesconsequences} for its precise relation with Bowen's original definition.
\item $X$ has the \emph{weak specification property} if there is $N\ge0$ such that for all $u,v\in\Lang(X)$ there is a word $w$ with $|w|\le N$ and $uwv\in\Lang(X)$. For collections of words, the analogous notion is the (W)-specification of \cite{climenhaga2012intrinsic}. It should not be confused with the weak specification property of \cite{kwietniak2015panoramaspecificationlikepropertiesconsequences}, where the connecting words may grow sublinearly with the length of the words being joined; $\Sigma(p_c)$ does not have that property either (Remark~\ref{rem:sublinear}).
\item $X$ is \emph{coded} if there is a countable collection $G$ of non-empty finite
words, called \emph{generators}, such that $X$ is the closure of the set of sequences
obtained by freely concatenating the generators. Since we work with one-sided
sequences, this means
\[
X=X_G:=\overline{\bigcup_{n\ge0}\sigma^n\big(\{g_1g_2g_3\cdots:\ g_i\in G\}\big)}.
\]
This is the one-sided version of the usual definition for two-sided shift spaces
\cite{blanchard1986systemes}; the shifts $\sigma^n$ are needed to make $X_G$
$\sigma$-invariant.
\item $X$ is \emph{synchronized} if it is transitive and there is $\omega\in\Lang(X)$ (a \emph{synchronizing word}) such that $u\omega,\ \omega v\in\Lang(X)$ imply $u\omega v\in\Lang(X)$. 
\item $X$ is \emph{hereditary} if $x\in X$ and $y\le x$ coordinatewise imply $y\in X$ \cite{kwietniak2011topological}.
\end{enumerate}
\end{definition}

Coded systems are transitive \cite{blanchard1986systemes}, and specification implies weak specification, which in turn implies transitivity. See \cite{kwietniak2015panoramaspecificationlikepropertiesconsequences} for a more detailed discussion of the relationships between different definitions of specification properties.
For $m,n\ge1$, every word of length $m+n$ in $\Lang(X)$ is the concatenation of a word of length $m$ and a word of length $n$, both in $\Lang(X)$, since the language of a subshift is closed under taking factors. Hence
\[
\#\Lang_{m+n}(X)\le\#\Lang_m(X)\cdot\#\Lang_n(X)\qquad(m,n\ge1),
\]
so the sequence $a_n=\log\#\Lang_n(X)$ is subadditive, and by Fekete's lemma the limit
\[
\lim_{n\to\infty}\frac1n\log\#\Lang_n(X)=\inf_{n\ge1}\frac1n\log\#\Lang_n(X)
\]
exists. This number is the \emph{topological entropy} of $X$ and is denoted by $h_{\mathrm{top}}(X)$. For an arbitrary collection of words $\mathcal D$, we write $h(\mathcal D)=\limsup_n\frac1n\log\#\mathcal D_n$, where $\mathcal D_n$ is the set of words of length $n$ in $\mathcal D$; since such a collection need not be closed under taking factors, a limit superior is used. 

For a $\sigma$-invariant Borel probability measure $\mu$ on $X$, we denote by $h_\mu(\sigma)$ its measure-theoretic entropy. The variational principle \cite[Thm.~8.6]{walters1982introduction} states that
\[
h_{\mathrm{top}}(X)=\sup\{h_\mu(\sigma):\ \mu \text{ is a $\sigma$-invariant Borel probability measure on } X\}.
\]
A \emph{measure of maximal entropy} (MME) is a $\sigma$-invariant Borel probability measure $\mu$ with $h_\mu(\sigma)=h_{\mathrm{top}}(X)$. Every subshift has at least one MME. Indeed, let $P$ be the partition of $X$ according to the first symbol, and for a finite partition $Q$, write 
\[
H_\mu(Q)=-\sum_{A\in Q}\mu(A)\log\mu(A). 
\]
The sets $\sigma^{-i}(A)$, $A\in P$, $i\ge0$, generate the Borel $\sigma$-algebra of $X$, so by the Kolmogorov--Sinai theorem for non-invertible maps \cite[Thm.~4.18]{walters1982introduction}, together with subadditivity,
\[
h_\mu(\sigma)=\inf_{n\ge1}\frac1n\,H_\mu\Bigl(\bigvee_{i=0}^{n-1}\sigma^{-i}P\Bigr).
\]
The partition $\bigvee_{i=0}^{n-1}\sigma^{-i}P$ consists of the cylinders of length $n$, which are clopen, so each term is a weak$^*$ continuous function of $\mu$, and $\mu\mapsto h_\mu(\sigma)$ is upper semicontinuous. Since the set of $\sigma$-invariant Borel probability measures on $X$ is weak$^*$ compact \cite[\S6.2]{walters1982introduction}, the supremum in the variational principle is attained. For expansive homeomorphisms, the upper semicontinuity of the entropy map is \cite[Thm.~8.2]{walters1982introduction}; the argument above is its one-sided version. The subshift $X$ is \emph{intrinsically ergodic} if it has exactly one MME. 

\subsection{Frequency shifts}

We recall the setting of \cite{becker2023symbolic}. For a non-empty word $w$, let $\delta_1(w)=\no{w}/|w|$.

\begin{definition}[\cite{becker2023symbolic}]
Let $p_c\in[0,1]$. A word $w=w_1\cdots w_n$ is \emph{admissible} if for every $1\le i<n$ with $w_i=1$, we have $w_{i+1}=1$ when $\delta_1(w_1\cdots w_i)\le p_c$, and $w_{i+1}=0$ when $\delta_1(w_1\cdots w_i)>p_c$. The set $P(p_c)$ consists of the sequences all of whose prefixes are admissible, and
\[
\Sigma(p_c)=\overline{\bigcup_{n\ge0}\sigma^n\big(P(p_c)\big)}.
\]
\end{definition}

The set $\Sigma(p_c)$ is a subshift: it is closed by definition, and it is $\sigma$-invariant because $\sigma$ is continuous and maps $\bigcup_{n\ge0}\sigma^n(P(p_c))$ into itself.

\subsection{The discrepancy}

The rule defining $P(p_c)$ compares the density $\delta_1(u)=\no{u}/|u|$ with $p_c$. To work with this condition, it is convenient to count the zeros and the ones of $u$ separately. Since $|u|=\nz{u}+\no{u}$, the density of a word $u$ with $\nz{u}$ zeros and $\no{u}$ ones is $\no{u}/(\nz{u}+\no{u})$, and this equals $p_c$ exactly when $\no{u}(1-p_c)=p_c\,\nz{u}$; if $\nz{u}>0$, this means that
\begin{equation*}
    \frac{\no{u}}{\nz{u}}=\frac{p_c}{1-p_c}.
\end{equation*}

The number
\begin{equation*}
    \beta=\frac{p_c}{1-p_c}\in(0,\infty)
\end{equation*}
is therefore the ratio between the numbers of ones and zeros of a word whose density is exactly $p_c$: at density $p_c$, each zero ``pays'' for $\beta$ ones. For instance, if $p_c=\frac25$, then $\beta=\frac23$, and a word with three zeros and two ones has density exactly $\frac25$. For $p_c=\frac{k}{k+1}$, we get $\beta=k$, that is, $k$ ones for each zero. 

This suggests measuring how far a word is from this balance by the difference between the number of ones that its zeros would allow and the number of ones it actually contains. For a word $w$, we define its \emph{discrepancy} by
\[
D(w)=\beta\,\nz{w}-\no{w}
\]
(in particular, $D(\varepsilon)=0$). Since $\nz{\cdot}$ and $\no{\cdot}$ are additive under concatenation, so is $D$: $D(uv)=D(u)+D(v)$ for all words $u,v$.
 
The role of $D$ in the transition rule is made precise in Lemmas~\ref{lem:density} and~\ref{lem:rule}: for a non-empty word $u$, the density satisfies $\delta_1(u)\le p_c$ if and only if $D(u)\ge0$, and consequently a $1$ is followed by $1$ if and only if the discrepancy of the prefix ending at it is non-negative.
\begin{lemma}[Discrepancy and density]\label{lem:density}
For every non-empty word $u$,
\[
(1-p_c)\,D(u)=|u|\,\bigl(p_c-\delta_1(u)\bigr).
\]
Consequently, $\delta_1(u)\le p_c$ if and only if $D(u)\ge0$, and $\delta_1(u)>p_c$ if and only if $D(u)<0$.
\end{lemma}

\begin{proof}
Since $\beta(1-p_c)=p_c$ and $|u|=\nz{u}+\no{u}$, we have
\[
(1-p_c)D(u)=p_c\,\nz{u}-(1-p_c)\,\no{u}=p_c\bigl(\nz{u}+\no{u}\bigr)-\no{u}=p_c|u|-\no{u}=|u|\bigl(p_c-\delta_1(u)\bigr).
\]
Both factors $1-p_c$ and $|u|$ are positive, so $D(u)$ and $p_c-\delta_1(u)$ are simultaneously non-negative or simultaneously negative.
\end{proof}

\begin{lemma}[Transition rule]\label{lem:rule}
A word $w$ is admissible if and only if, for every index $i$ with $1\le i<|w|$ and $w_i=1$,
\begin{itemize}
\item $w_{i+1}=1$ whenever $D(w_1\cdots w_i)\ge0$, and
\item $w_{i+1}=0$ whenever $D(w_1\cdots w_i)<0$.
\end{itemize}
In particular, for such $i$ we have $w_{i+1}=1$ if and only if $D(w_1\cdots w_i)\ge0$.
\end{lemma}

\begin{proof}
Fix $i$ with $1\le i<|w|$ and $w_i=1$. The prefix $w_1\cdots w_i$ is non-empty, so by Lemma~\ref{lem:density} the hypothesis $\delta_1(w_1\cdots w_i)\le p_c$ of the definition of admissibility is the hypothesis $D(w_1\cdots w_i)\ge0$, and the hypothesis $\delta_1(w_1\cdots w_i)>p_c$ is the hypothesis $D(w_1\cdots w_i)<0$. Hence the condition imposed on $w_{i+1}$ in the definition of admissibility at position $i$ coincides with the two conditions in the statement. The last assertion holds because exactly one of $D(w_1\cdots w_i)\ge0$ and $D(w_1\cdots w_i)<0$ is true and $w_{i+1}\in\{0,1\}$: if $D\ge0$, then $w_{i+1}=1$, and if $D<0$, then $w_{i+1}=0\ne1$.
\end{proof}

Thus each $0$ increases the discrepancy by $\beta$, each $1$ decreases it by $1$, and a run of ones continues as long as the discrepancy stays non-negative. The next lemma makes this precise: the length of a run of ones is determined by the discrepancy at its start. This is the basic computation of the paper. It leads to the estimate on anchored runs (Lemma~\ref{lem:anchored}), and through it to Theorems~\ref{thm:A} and~\ref{thm:B}; it also gives the structure of $P_k$ (Lemma~\ref{lem:structure}). 

\begin{lemma}[Run lengths]\label{lem:runs}
Let $w$ be admissible, let $w_{[i,j]}=1^r$ be a maximal run of ones in $w$, and put $d=D(w_1\cdots w_{i-1})$.
\begin{enumerate}[label=\textup{(\alph*)}]
\item If the run is closed, then $r=\lfloor d\rfloor+1$ when $d\ge0$, and $r=1$ when $d<0$. Moreover, if $d\ge0$, then $D(w_1\cdots w_j)=d-\lfloor d\rfloor-1\in[-1,0)$.
\item If the run is not closed (that is, $j=|w|$), then $r\le\lfloor d\rfloor+1$ when $d\ge0$, and $r=1$ when $d<0$.
\end{enumerate}
\end{lemma}

\begin{proof}
Each $1$ contributes $-1$ to $D$, so, since $w_i,\dots,w_{i+t-1}$ are ones, additivity gives
\begin{equation}\label{eq:dt}
D(w_1\cdots w_{i+t-1})=d-t\qquad(1\le t\le r).
\end{equation}
Note that $w_{i+t-1}$ is the $t$-th one of the run. Two facts follow from Lemma~\ref{lem:rule}.
\begin{enumerate}[label=(\roman*)]
\item If $1\le t<r$, then $w_{i+t-1}=1$ and $w_{i+t}=1$, so the rule forces $D(w_1\cdots w_{i+t-1})\ge0$; by \eqref{eq:dt}, $t\le d$.
\item If the run is closed, then $w_j=1$ and $w_{j+1}=0$ with $j=i+r-1<|w|$, so the rule forces $D(w_1\cdots w_j)<0$; by \eqref{eq:dt}, $d<r$.
\end{enumerate}
\emph{Case $d<0$.} If $r\ge2$, then (i) with $t=1$ gives $1\le d$, which is absurd. Hence $r=1$, in both (a) and (b).

\emph{Case $d\ge0$.} If $r\ge2$, then (i) with $t=r-1$ gives $r-1\le d$; if $r=1$ this holds trivially. Since $r-1$ is an integer, $r\le\lfloor d\rfloor+1$ in all cases, which proves (b). If the run is closed, (ii) adds $d<r$, so $r-1\le d<r$, that is, $r-1=\lfloor d\rfloor$. Finally, \eqref{eq:dt} with $t=r$ gives $D(w_1\cdots w_j)=d-r=d-\lfloor d\rfloor-1$, and $r-1\le d<r$ shows that this lies in $[-1,0)$. This proves (a).
\end{proof}

\begin{example}\label{ex:discrepancy}
Let $p_c=\frac25$, so $\beta=\frac23$. Choosing zero blocks of lengths $3,2,1,4$ forces the admissible word
\[
w=0^3\,1^3\;0^2\,1\;0\,1\;0^4\,1^2\;0,\qquad |w|=18.
\]
Its maximal runs of ones are $w_{[4,6]}=111$, $w_{[9,9]}=1$, $w_{[11,11]}=1$ and $w_{[16,17]}=11$, and all of them are closed; in the prefix $w_{[1,17]}$, the last run $w_{[16,17]}$ is not closed. The factor $w_{[4,5]}=11$ consists of ones but is not a maximal run, because $w_6=1$. The discrepancies at the start of the runs of ones are $2,\ \frac13,\ 0,\ \frac53$, so the runs have lengths $\lfloor d\rfloor+1=3,1,1,2$ (Lemma~\ref{lem:runs}(a)), and the discrepancies at their ends are $-1,-\frac23,-1,-\frac13$.
\end{example}

\begin{lemma}[Language]\label{lem:language}
\begin{enumerate}[label=\textup{(\alph*)}]
\item Every admissible word is a prefix of some $x\in P(p_c)$.
\item $\Lang(\Sigma(p_c))$ is the set of factors of admissible words, equivalently the set of factors of elements of $P(p_c)$.
\end{enumerate}
\end{lemma}

\begin{proof}
(a) Let $w$ be admissible. If $w$ is empty or ends with $0$, or if $w$ ends with $1$ and $D(w)<0$, put $s=0$. Otherwise, $w$ ends with $1$ and $e:=D(w)\ge0$, and we put $s=\lfloor e\rfloor+1$. We claim that $x=w1^s0^\infty$ belongs to $P(p_c)$.

A prefix of $x$ of length $n$ is admissible if and only if the condition of Lemma~\ref{lem:rule} holds at every $i<n$ with $x_i=1$; this condition involves only $x_1\cdots x_{i+1}$. So it suffices to check it at every $i$ with $x_i=1$. For $i<|w|$, it holds because $w$ is admissible. Positions $i>|w|+s$ carry zeros. It remains to consider the positions $i=|w|+t$ with $0\le t\le s$ in the case $s\ge1$, where $x_i=1$ and $D(x_1\cdots x_i)=e-t$. For $t<s$, we have $t\le\lfloor e\rfloor$, so $e-t\ge0$ and $x_{i+1}=1$, as required. For $t=s$, we have $e-s=e-\lfloor e\rfloor-1<0$ and $x_{i+1}=0$, as required. (If $s=0$ and $w$ ends with $1$, then $D(w)<0$ and $x_{|w|+1}=0$, as required.)

(b) Elements of $P(p_c)$ lie in $\Sigma(p_c)$, so their factors belong to $\Lang(\Sigma(p_c))$. Every admissible word is a prefix of an element of $P(p_c)$ by (a), and every prefix of an element of $P(p_c)$ is admissible; hence the factors of admissible words are exactly the factors of elements of $P(p_c)$. Conversely, let $u$ occur at position $i$ in some $y\in\Sigma(p_c)$. The set of sequences with $u$ at position $i$ is open, so, $y$ being in the closure, it contains a point $\sigma^n(x)$ with $x\in P(p_c)$. Then $u$ is a factor of $x$.
\end{proof}

The next lemma is the key estimate. A single $1$ followed by $0$ only gives an upper bound on the discrepancy, but a pair $11$ followed by $0$ pins it down to an interval of length one. 
\begin{lemma}[Anchored runs]\label{lem:anchored}
Let $m,n\ge1$ with $\beta m\ge1$. If $110^m1^n0\in\Lang(\Sigma(p_c))$, then
\[
\beta m-1<n<\beta m+1.
\]
\end{lemma}

\begin{proof}
By Lemma~\ref{lem:language}(b), there is an admissible word $w=a\,11\,0^m1^n0\,b$. Put $D_0=D(a11)$. The first $1$ of the pair $11$ is followed by a $1$, so $D(a1)\ge0$ by Lemma~\ref{lem:rule}, and $D_0=D(a1)-1\ge-1$. The second $1$ is followed by a $0$, so $D_0<0$ by Lemma~\ref{lem:rule}. Since $m\ge1$ and the block $1^n$ is followed by $0$, the block $1^n$ is a closed maximal run of ones; it starts after the prefix $a110^m$, whose discrepancy is $d=D_0+\beta m\in[\beta m-1,\beta m)$. As $\beta m\ge1$, we have $d\ge0$, and Lemma~\ref{lem:runs}(a) gives $n=\lfloor d\rfloor+1$, hence $d<n\le d+1$. Therefore $\beta m-1\le d<n\le d+1<\beta m+1$.
\end{proof}

\begin{lemma}\label{lem:existence}
For every $m\ge1$, there exist $a,r,n\ge1$ with $r\ge2$ such that $0^a1^r0^m1^n0$ is admissible. In particular, there is $n_m\ge1$ such that
\[ 110^m1^{n_m}0\in\Lang(\Sigma(p_c)),\]
and hence $110^m\in\Lang(\Sigma(p_c))$.
\end{lemma}

\begin{proof}
Let $a=\lceil1/\beta\rceil$, so that $\beta a\ge1$, and let $r=\lfloor\beta a\rfloor+1\ge2$. The word $0^a1^r0$ is admissible
by Lemma~\ref{lem:rule}: the discrepancy of the prefix ending at the $t$-th one is $D(0^a1^t)=\beta a-t$, which is $\ge0$ for $t<r$ (as $t\le\lfloor\beta a\rfloor$), while $D(0^a1^r)=\beta a-r<0$.
Next, $w'=0^a1^r0^m1$ is admissible, because $m\ge1$ implies that its last symbol is not read by any condition of Lemma~\ref{lem:rule}, so $w'$ satisfies the same conditions as $0^a1^r0^{m+1}$. By the proof of Lemma~\ref{lem:language}(a), $w'$ is a prefix of $x=w'1^s0^\infty\in P(p_c)$ for some $s\ge0$. Put $n=s+1$. Then $0^a1^r0^m1^n0$ is a prefix of $x$, hence admissible. Since $r\ge2$, the word $110^m1^n0$ is a factor of it (the last two ones of $1^r$, followed by $0^m1^n0$). By Lemma~\ref{lem:language}(b), $110^m1^n0\in\Lang(\Sigma(p_c))$, and so is its prefix $110^m$.
\end{proof}

In Lemma~\ref{lem:anchored}, the length $n$ is \emph{not} a function of $m$, because it depends on the discrepancy $D_0=D(a11)\in[-1,0)$ at the anchoring pair (see the proof). The following example shows this.

\begin{example}\label{ex:anchored}
Let $p_c=\frac13$, so that $\beta=\frac12$, and let $m\ge2$. The words $0^21^2$ and $0^31^2$ are admissible and end with discrepancy $D_0=-1$ and $D_0=-\frac12$, respectively. After a further block $0^m$, the discrepancy becomes $d=\frac m2-1$ in the first case and $d=\frac{m-1}2$ in the second, and the next run of ones has length $\lfloor d\rfloor+1$. Checking the rule of Lemma~\ref{lem:rule} at each $1$, one sees that the words
\[
0^21^2\,0^m1^{\lfloor m/2\rfloor}0
\qquad\text{and}\qquad
0^31^2\,0^m1^{\lfloor (m+1)/2\rfloor}0
\]
are admissible. Hence $110^m1^n0\in\Lang(\Sigma(\frac13))$ both for $n=\lfloor\frac m2\rfloor$ and for $n=\lfloor\frac{m+1}2\rfloor$. These values coincide when $m$ is even and differ when $m$ is odd; for instance, both $110^3\,1\,0$ and $110^3\,1^2\,0$ belong to $\Lang(\Sigma(\frac13))$. By Lemma~\ref{lem:anchored}, no other value of $n$ is possible. This is why the proof of Theorem~\ref{thm:A} compares lengths $m<m'$ with $\beta(m'-m)\ge2$.
\end{example}

\begin{proposition}\label{prop:not-hereditary}
For every $p_c\in(0,1)$, the subshift $\Sigma(p_c)$ is not hereditary.
\end{proposition}

\begin{proof}
Suppose, for a contradiction, that $\Sigma(p_c)$ is hereditary. Let $m\ge1$ with $\beta m\ge3$, and let $n=n_m$ be given by Lemma~\ref{lem:existence}, so that $110^m1^n0\in\Lang(\Sigma(p_c))$. By Lemma~\ref{lem:anchored}, $\beta m-1<n<\beta m+1$; in particular $n>\beta m-1\ge2$, so $n\ge3$.

Let $x\in\Sigma(p_c)$ be a point in which $110^m1^n0$ occurs, say $x_ix_{i+1}\cdots x_{i+m+n+2}=110^m1^n0$, and let $y$ be obtained from $x$ by replacing the last two symbols of this block $1^n$ by zeros; that is, $y_j=0$ for $j\in\{i+m+n,\,i+m+n+1\}$ and $y_j=x_j$ otherwise. Then $y\le x$ coordinatewise, so $y\in\Sigma(p_c)$ by our assumption. Hence every word occurring in $y$ belongs to $\Lang(\Sigma(p_c))$; in particular
$110^m1^{n-2}0\in\Lang(\Sigma(p_c))$. Since $n-2\ge1$, Lemma~\ref{lem:anchored} applies to this word (with $n-2$ in place of $n$) and yields $\beta m-1<n-2$, that is, $n>\beta m+1$. This contradicts $n<\beta m+1$. Therefore $\Sigma(p_c)$ is not hereditary.
\end{proof}

%==================================================================
\section{Non-soficity}\label{sec:nonsofic}
%==================================================================

Sofic subshifts were introduced by Weiss \cite{weiss1973subshifts}. Most of
the literature, including the standard reference \cite{Marcus}, works with
two-sided subshifts, that is, with spaces of bi-infinite sequences. Since
$\Sigma(p_c)$ is a one-sided subshift, we use the definition given by Ashley,
Kitchens and Stafford \cite[\S2]{ashley1992boundaries}, who explicitly treat
one-sided subshifts and characterize soficity in terms of follower sets.

\begin{definition}[{\cite[\S2]{ashley1992boundaries}}]\label{def:sofic}
A subshift $X\subseteq\{0,1\}^{\N}$ is \emph{sofic} if the collection
$\{\Fol_X(w):w\in\Lang(X)\}$ is finite.
\end{definition}

As noted in \cite[\S2]{ashley1992boundaries}, this is equivalent to $X$ being
a factor of a one-sided subshift of finite type. 
\begin{theorem}[Theorem~\ref{thm:A}]\label{thm:nonsofic}
Let $p_c\in(0,1)$ and $m_i=i\lceil2/\beta\rceil$ for $i\ge1$. The follower sets
$\Fol(110^{m_i})$ in $\Sigma(p_c)$ are pairwise distinct. Consequently,
$\Sigma(p_c)$ is not sofic.
\end{theorem}

\begin{proof}
Note that $\beta m_i\ge2$ for all $i$, and $\beta(m_j-m_i)\ge2$ for $i<j$. Fix
$i<j$ and let $n=n_{m_i}$ be given by Lemma~\ref{lem:existence}, so that
$1^n0\in\Fol(110^{m_i})$. By Lemma~\ref{lem:anchored}, $n<\beta m_i+1$. If
$1^n0\in\Fol(110^{m_j})$, then Lemma~\ref{lem:anchored} would give
$n>\beta m_j-1\ge\beta m_i+1$, a contradiction. Hence
$\Fol(110^{m_i})\ne\Fol(110^{m_j})$, so $\Sigma(p_c)$ has infinitely many
distinct follower sets and is not sofic by Definition~\ref{def:sofic}.
\end{proof}

\begin{corollary}[cf.~{\cite[Thm.~1]{becker2023symbolic}}]\label{cor:notSFT}
For every $p_c\in(0,1)$, $\Sigma(p_c)$ is not a subshift of finite type.
\end{corollary}
\begin{proof}
In a subshift of finite type whose forbidden words have length at most $N$, the follower set of a word of length at least $N-1$ depends only on its last $N-1$ symbols, so there are only finitely many follower sets. Hence every subshift of finite type is sofic, and the corollary follows from Theorem~\ref{thm:nonsofic}.
\end{proof}

\begin{remark}\label{rem:BB-thm1}
The proof of \cite[Thm.~1]{becker2023symbolic} uses the fact that blocks $0^m1^n0$ with $n-\beta m<-1$ do not occur in $\Sigma(p_c)$. This is only true under an additional assumption. The left-hand inequalities of \cite[(3), (5)]{becker2023symbolic} express the admissibility condition at the position preceding the last $1$ of the block placed before $0^m1^n0$, and that position carries a transition $1\to1$ only when the corresponding run of ones has length at least $2$. Without it, the conclusion may fail: for $p_c=\frac13$ ($\beta=\frac12$), the word $(10)^6\,0^5\,1\,0$ is admissible, so $0^510\in\Lang(\Sigma(\frac13))$ although $n-\beta m=-\frac32$. Anchoring with $11$, as in Lemma~\ref{lem:anchored}, supplies exactly the missing hypothesis, and Corollary~\ref{cor:notSFT} gives a complete proof of \cite[Thm.~1]{becker2023symbolic}.
\end{remark}

%==================================================================
\section{Absence of specification}\label{sec:nospec}
%==================================================================

\begin{theorem}[Theorem~\ref{thm:B}]\label{thm:nospec}
Let $p_c\in(0,1)$ and let $n\ge1$ with $\beta n\ge1$. If $w$ is a word such that $110^n\,w\,10\in\Lang(\Sigma(p_c))$, then $|w|>\beta n-2$. In particular, $\Sigma(p_c)$ does not have the weak specification property.
\end{theorem}

\begin{proof}
Both $u=110^n$ and $v=10$ belong to $\Lang(\Sigma(p_c))$: the first by Lemma~\ref{lem:existence}, and the second because $D(1)=-1<0$. Suppose $uwv$ is a factor of an admissible word $y$ (Lemma~\ref{lem:language}). In $y$, the maximal run of zeros that follows the pair $11$ of $u$ has length $m\ge n$ and ends before the $1$ of $v$. Let $R$ be the maximal run of ones that follows it. Then $R$ starts at or before the $1$ of $v$ and ends at or before it; in particular, $R$ is closed. Hence $R$ is contained in the segment occupied by $w1$, and $110^m1^{|R|}0$ is a factor of $y$. Since $\beta m\ge\beta n\ge1$, Lemma~\ref{lem:anchored} gives $|R|>\beta m-1\ge\beta n-1$. Therefore $|w|+1\ge|R|>\beta n-1$.

Suppose that $\Sigma(p_c)$ has the weak specification property with constant $N$. Since $\beta>0$, we can choose $n\ge1$ with $\beta n\ge N+2$; in particular, $\beta n\ge1$. Applying weak specification to $u=110^n$ and $v=10$, we obtain a word $w$ with $|w|\le N$ and $110^n\,w\,10\in\Lang(\Sigma(p_c))$. By the first part of the proof, $|w|>\beta n-2\ge N$, a contradiction.
\end{proof}

Since specification implies weak specification, $\Sigma(p_c)$ does not have the specification property either.

\begin{remark}\label{rem:sublinear}
Theorem~\ref{thm:nospec} also rules out the weak specification property of \cite{kwietniak2015panoramaspecificationlikepropertiesconsequences}, in which the connecting words may grow sublinearly with the length of the words being joined. Indeed, for every $n$ with $\beta n\ge1$, joining $110^n$ to $10$ requires a connecting word of length larger than $\beta n-2$. This bound tends to infinity although $10$ has fixed length, and it is linear in the length $n+2$ of $110^n$; hence no sublinear gap function can work.
\end{remark}

Although our standing assumption is $p_c\in(0,1)$, it is instructive to compare with the endpoints. For $p_c=0$, we have $\beta=0$, so the condition $\beta n\ge1$ of Theorem~\ref{thm:nospec} never holds, consistent with $\Sigma(0)$ being the golden mean shift, which has specification. For $p_c=1$, the parameter $\beta$ is not defined, and the obstruction is of a different nature: $\Sigma(1)=\{0^\infty\}\cup\{0^j1^\infty:j\ge0\}$ is not even transitive.

%==================================================================
\section{The parameters \texorpdfstring{$p_c=\frac{k}{k+1}$}{pc=k/(k+1)}}\label{sec:k}
%==================================================================

Throughout this section, $k\ge1$ and $p_c=\frac{k}{k+1}$, so that $\beta=k$, and we write $\Sigma_k=\Sigma(\frac{k}{k+1})$ and $P_k=P(\frac{k}{k+1})$. The discrepancy is now integer-valued.

\begin{lemma}[Structure of $P_k$; cf.~{\cite[Cor.~1, 2]{becker2023symbolic}}]\label{lem:structure}
A sequence $x$ belongs to $P_k$ if and only if it has one of the forms
\begin{enumerate}[label=\textup{(\roman*)}]
\item $x=0^\infty$;
\item $x=0^{m_1}1^{km_1+1}\,0^{m_2}1^{km_2}\,0^{m_3}1^{km_3}\cdots$;
\item $x=1\,0^{m_1}1^{km_1}\,0^{m_2}1^{km_2}\cdots$;
\end{enumerate}
with all $m_i\ge1$, where in (ii) and (iii) the sequence of blocks may be finite (in (iii), possibly empty: $x=10^\infty$) and followed by $0^\infty$. Moreover, the discrepancy equals $-1$ at the end of every closed run of ones.
\end{lemma}

\begin{proof}
Every run of ones of $x$ is finite, by Lemma~\ref{lem:runs}(b) applied to the prefixes of $x$; hence every run is closed, and Lemma~\ref{lem:runs}(a) applies to it. If $x$ starts with $0^{m_1}1$, the first run starts with $d=km_1\ge1$, so it has length $km_1+1$ and ends with discrepancy $-1$. If $x_1=1$, then $D(1)=-1<0$ and the run has length $1$. If a run ends with discrepancy $-1$ and is followed by $0^m1$ with $m\ge1$, the next run starts with $d=km-1\ge0$, so it has length $km$ and ends again with discrepancy $-1$. Conversely, sequences of these forms satisfy Lemma~\ref{lem:rule}, since zeros are unconstrained and every run has the forced length.
\end{proof}

\begin{definition}
Let $G_k=\{0^m1^{km}:m\ge1\}$, let $G_k^*$ be the set of finite concatenations of elements of $G_k$ (including $\eps$), let $\Lang_k$ be the set of factors of words in $G_k^*$, and let $X_k=X_{G_k}$. A \emph{bad word} is a word $01^r0$ with $k\nmid r$.
\end{definition}

When $k=1$, there are no bad words. We have $\Lang(X_k)=\Lang_k$. Indeed, every word of $\Lang_k$ is a factor of a sequence $g_1g_2g_3\cdots$ with $g_i\in G_k$, and such a sequence lies in $X_k$. Conversely, if a word occurs in a point of $X_k$, then, by the argument in the proof of Lemma~\ref{lem:language}(b), it occurs in some sequence $g_1g_2g_3\cdots$, and hence in some word $g_1\cdots g_\ell\in G_k^*$. 
\begin{lemma}\label{lem:boundaries}
Let $w=g_1\cdots g_\ell\in G_k^*$ with $g_i=0^{m_i}1^{km_i}\in G_k$. The maximal runs of $w$ are, from left to right, $0^{m_1},1^{km_1},\dots,0^{m_\ell},1^{km_\ell}$. In particular, the occurrences of $10$ in $w$ are exactly the \emph{junctions}, formed by the last letter of $g_i$ and the first letter of $g_{i+1}$ $(1\le i<\ell)$, and no word of $\Lang_k$ contains a bad word.
\end{lemma}

\begin{proof}
Since $m_i\ge1$, each $g_i$ is a non-empty block of zeros followed by a non-empty block of ones, so the blocks listed in the statement alternate between zeros and ones. This gives the first claim, and an occurrence of $10$ is then the end of some $1^{km_i}$ followed by the beginning of $0^{m_{i+1}}$. If a bad word $01^r0$ occurred in a word of $\Lang_k$, it would occur in some $w\in G_k^*$, and $1^r$ would be a maximal run of ones of $w$; then $r=km_i$ for some $i$, contradicting $k\nmid r$.
\end{proof}

\begin{proposition}\label{prop:language-k}
\begin{enumerate}[label=\textup{(\alph*)}]
\item $\Lang_k\subseteq\Lang(\Sigma_k)$; hence $X_k\subseteq\Sigma_k$.
\item If $u\in\Lang(\Sigma_k)$ contains no bad word, then $u\in\Lang_k$.
\item A point $x\in\Sigma_k$ belongs to $X_k$ if and only if it contains no bad word.
\end{enumerate}
\end{proposition}

\begin{proof}

(a) Let $g_1\cdots g_\ell\in G_k^*$ with $g_i=0^{m_i}1^{km_i}$. The sequence $1\,g_1\cdots g_\ell\,0^\infty$ is of form (iii) in Lemma~\ref{lem:structure} (the leading $1$ makes the initial discrepancy equal to $-1$, so that the first run has length $km_1$), hence it belongs to $P_k$ and all its prefixes are admissible; in particular $1g_1\cdots g_\ell$ is admissible. By Lemma~\ref{lem:language}(b), $g_1\cdots g_\ell\in\Lang(\Sigma_k)$. Since $\Lang(\Sigma_k)$ is closed under taking factors, $\Lang_k\subseteq\Lang(\Sigma_k)$. Finally, every factor of a point of $X_k$ lies in $\Lang_k\subseteq\Lang(\Sigma_k)$, and a sequence all of whose factors lie in $\Lang(\Sigma_k)$ belongs to $\Sigma_k$; hence $X_k\subseteq\Sigma_k$.

(b) Let $u$ be a factor of $x\in P_k$ (Lemma~\ref{lem:language}). If $x=0^\infty$, then $u=0^a$ is a prefix of $0^a1^{ka}$. If $x$ is of form (iii), then $x$ is a suffix of $01^k\,0^{m_1}1^{km_1}\cdots$, so every factor of $x$ lies in $\Lang_k$ (a final $0^a$ is a prefix of $0^a1^{ka}$). If $x=0^{m_1}1^{km_1+1}y$ is of form (ii), then:
\begin{itemize}
\item every factor of $0^{m_1}1^{km_1+1}$ is a factor of $0^{m_1+1}1^{k(m_1+1)}$;
\item $1^{km_1+1}y$ is a suffix of $0^{m_1+1}1^{k(m_1+1)}y$, because $km_1+1\le k(m_1+1)$, so all its factors lie in $\Lang_k$.
\end{itemize}
The only remaining factors of $x$ are those containing $0\,1^{km_1+1}0$, the first run together with both neighbors. For $k\ge2$\corr{,} this is a bad word. For $k=1$, $x$ is a suffix of $0x=0^{m_1+1}1^{m_1+1}y$, and all its factors lie in $\Lang_1$.

(c) If $x\in\Sigma_k$ has no bad word, all its factors lie in $\Lang_k$ by (b), so $x\in X_k$. The converse follows from Lemma~\ref{lem:boundaries}.
\end{proof}

\begin{proposition}\label{prop:Xk-properties}
For every $k\ge1$, $X_k$ is coded, transitive and synchronized, with synchronizing word $10$.
\end{proposition}
\begin{proof}
$X_k$ is coded by definition. For transitivity, let $u,v\in \Lang_k$ . Choose $a,b\in G_k^*$ and words $a',a'',b',b''$ with $a=a'ua''$ and $b=b'vb''$. Since $G_k^*$ is closed under concatenation, $ab=a'\,u\,(a''b')\,v\,b''\in G_k^*,$ so $uwv\in\Lang_k$ with $w=a''b'$, and $X_k$ is transitive. Now let $u10,\,10v\in\Lang_k$. By Lemma~\ref{lem:boundaries}, $u10$ is a factor of $a_1a_2$ with $a_1,a_2\in G_k^*$, where $a_1$ ends with the $1$ and $a_2$ starts with the $0$ of this occurrence of $10$. Thus $u1$ is a suffix of $a_1$. Likewise $0v$ is a prefix of some $b_2\in G_k^*$. Hence $u10v$ is a factor of $a_1b_2\in G_k^*$.
\end{proof}

\begin{remark}\label{rem:period}
Every generator $0^m1^{km}$ has length divisible by $k+1$, and by Lemma~\ref{lem:boundaries} the occurrences of $10$ in a word of $G_k^*$ are exactly the junctions between generators. Hence, if $10\,w\,10\in\Lang_k$, then $|w|+2$ is a sum of lengths of generators, so $|w|\equiv-2\pmod{k+1}$. In particular, $X_k$ is not topologically mixing for any $k\ge1$; (recall that a subshift $X$ is \emph{topologically mixing} if for all $u,v\in\Lang(X)$ there is $N$ such that for every $n\ge N$ there is a word $w$ with $|w|=n$ and $uwv\in\Lang(X)$) for instance, $\Sigma(\frac12)=X_1$ is transitive and synchronized, but $10\,w\,10\in\Lang(\Sigma(\frac12))$ forces $|w|$ to be even.

\end{remark}

\begin{example}[Theorem~\ref{thm:C}(a)]\label{ex:half}
For $p_c=\frac12$, we have $\beta=1$, and $\Sigma(\frac12)=X_1$ is the coded system generated by $\{0^m1^m:m\ge1\}$. Indeed, for $k=1$, there are no bad words, so Proposition~\ref{prop:language-k} gives $\Lang(\Sigma(\tfrac12))=\Lang_1$. By Proposition~\ref{prop:Xk-properties}, $\Sigma(\frac12)$ is therefore coded, transitive and synchronized, with synchronizing word $10$.
\end{example}

\begin{proposition}[Theorem~\ref{thm:C}(b)]\label{prop:k2}
Let $k\ge2$.
\begin{enumerate}[label=\textup{(\alph*)}]
\item Every word of $\Lang(\Sigma_k)$, and hence every point of $\Sigma_k$, contains at most one occurrence of a bad word.
\item $\Sigma_k$ is not transitive. In particular, it is neither coded nor synchronized.
\item $X_k\subsetneq\Sigma_k$.
\end{enumerate}
\end{proposition}

\begin{proof}
(a) An occurrence of a bad word in a factor of $x\in P_k$ is a maximal run of ones in $x$ preceded and followed by $0$, with length not divisible by $k$. By Lemma~\ref{lem:structure}, the only such run is the first run of a sequence of form (ii), whose length $km_1+1\equiv1\pmod k$. In form (iii), the first run sits at position $1$, so it is not preceded by $0$. If a point had two occurrences, a common factor containing both would lie in $\Lang(\Sigma_k)$, contradicting what we have just shown.

(b) The words $u=01^{k+1}0$ and $v=0^21^{2k+1}0$ are prefixes of sequences of form (ii), with $m_1=1,2$, so they lie in $\Lang(\Sigma_k)$. Both are bad, since $k+1\equiv2k+1\equiv1\pmod k$. By (a), no word $uwv$ lies in $\Lang(\Sigma_k)$.

(c) $01^{k+1}0\in\Lang(\Sigma_k)\setminus\Lang_k$.
\end{proof}

\begin{proposition}\label{prop:complement}
Let $k\ge2$, and let $Y_k$ be the set of sequences $g_1g_2g_3\cdots$ and $g_1\cdots g_\ell\,0^\infty$ $(\ell\ge0)$ with $g_i\in G_k$. Then $Y_k$ is a closed subset of $X_k$, and
\[
\Sigma_k\setminus X_k=\bigl\{0^j1^{km+1}y:\ m\ge1,\ 1\le j\le m,\ y\in Y_k\bigr\}.
\]
In particular, if $a\,01^{km+1}0\in\Lang(\Sigma_k)$ for some word $a$ and some $m\ge1$, then $a=0^i$ with $0\le i\le m-1$; for instance, $01^{k+1}0\in\Lang(\Sigma_k)$, but neither $001^{k+1}0$ nor $101^{k+1}0$ belongs to $\Lang(\Sigma_k)$; that is, $01^{k+1}0$ can only occur at the first coordinate of a point of $\Sigma_k$.
\end{proposition}

\begin{proof}
A sequence belongs to $Y_k$ if and only if it begins with $0$ and every finite maximal run of zeros $0^c$ in it is followed by $1^{kc}0$. Each of these conditions involves finitely many coordinates, so $Y_k$ is closed. The sequences $g_1g_2g_3\cdots$ lie in $X_k$ by definition, and $g_1\cdots g_\ell\,0^\infty$ is the limit of $g_1\cdots g_\ell\,(0^n1^{kn})^\infty$ as $n\to\infty$; hence $Y_k\subseteq X_k$.

The sequences of form~(ii) in Lemma~\ref{lem:structure} are exactly the sequences $0^m1^{km+1}y$ with $m\ge1$ and $y\in Y_k$. Hence the points on the right-hand side, $0^j1^{km+1}y=\sigma^{m-j}(0^m1^{km+1}y)$, lie in $\Sigma_k$. They contain the bad word $01^{km+1}0$, because $j\ge1$, $y$ begins with $0$, and $k\nmid km+1$ (as $k\ge2$). By Proposition~\ref{prop:language-k}(c), they do not lie in $X_k$.

Conversely, let $x\in\Sigma_k\setminus X_k$. By Proposition~\ref{prop:language-k}(c), $x$ contains a bad word, say $x_{[j,j+r+1]}=01^r0$ with $k\nmid r$. Let $N\ge j+r+1$. By Lemma~\ref{lem:language}(b), $x_{[1,N]}$ is a factor of some $z\in P_k$. By Lemma~\ref{lem:structure}, the only maximal run of ones of a point of $P_k$ that is preceded and followed by $0$ and has length not divisible by $k$ is the first run of a sequence of form~(ii) (see the proof of Proposition~\ref{prop:k2}(a)). Hence $z=0^m1^{km+1}y'$ with $m\ge1$ and $y'\in Y_k$, the run $1^r$ of $x_{[1,N]}$ occupies the first run of $z$, and $r=km+1$; in particular, $m$ does not depend on $N$. Since the first run of $z$ starts at coordinate $m+1$, we get $x_{[1,N]}=z_{[m-j+1,\,m-j+N]}$. Therefore $j\le m$, $x_{[1,j]}=0^j$, and $x_{[j+r+1,N]}$ is a prefix of $y'$. Since $N$ is arbitrary and $Y_k$ is closed, $\sigma^{j+r}(x)\in Y_k$, and so $x=0^j1^{km+1}\sigma^{j+r}(x)$ belongs to the right-hand side.

Finally, if $a\,01^{km+1}0\in\Lang(\Sigma_k)$, apply the previous paragraph to a point $x\in\Sigma_k$ beginning with this word, with $j=|a|+1$ and $r=km+1$. This gives $|a|+1\le m$ and $a=x_{[1,|a|]}=0^{|a|}$.
\end{proof}

In particular, if $x=0^j1^{km+1}y\in\Sigma_k\setminus X_k$ as above, then $\sigma^j(x)=1^{km+1}y$ contains no bad word, so $\sigma^j(x)\in X_k$ by Proposition~\ref{prop:language-k}(c): every point of $\Sigma_k\setminus X_k$ enters $X_k$ after finitely many steps. Since $\sigma(X_k)\subseteq X_k$, the sets $W_j=\sigma^{-j}(X_k)\cap\Sigma_k$ increase with $j$ and cover $\Sigma_k$, and $\mu(W_j)=\mu(X_k)$ for every $\sigma$-invariant measure $\mu$; hence $\mu(X_k)=\mu(\Sigma_k)=1$. This gives a second proof of Proposition~\ref{prop:measures}.

\begin{proposition}[Theorem~\ref{thm:C}(b), measures]\label{prop:measures}
For every $k\ge1$, every $\sigma$-invariant Borel probability measure $\mu$ on $\Sigma_k$ satisfies $\mu(X_k)=1$. Consequently, $h_{\mathrm{top}}(\Sigma_k)=h_{\mathrm{top}}(X_k)$, and $\Sigma_k$ and $X_k$ have the same measures of maximal entropy.
\end{proposition}

\begin{proof}
For $k=1$, there is nothing to prove, since $\Sigma_1=X_1$ (Example~\ref{ex:half}). Let $k\ge2$, and let $A_n$ be the set of $x\in\Sigma_k$ having a bad word that starts at coordinate $n$. Each $A_n$ is a union of cylinders, hence Borel, and $\sigma^{-1}(A_n)\cap\Sigma_k=A_{n+1}$. By Proposition~\ref{prop:k2}(a), the sets $A_n$ are pairwise disjoint, while invariance gives $\mu(A_n)=\mu(A_1)$ for all $n$. Hence $\mu(A_1)=0$, so $\mu\big(\bigcup_nA_n\big)=0$, and $\Sigma_k\setminus\bigcup_nA_n=X_k$ by Proposition~\ref{prop:language-k}(c). The remaining statements follow from the variational principle \cite[Thm.~8.6]{walters1982introduction}.
\end{proof}

%==================================================================
\section{Intrinsic ergodicity for \texorpdfstring{$p_c=\frac{k}{k+1}$}{pc=k/(k+1)}}\label{sec:IE}
%==================================================================
\subsection{The Climenhaga--Thompson criterion}

Let $X$ be a subshift with language $\Lang$. A collection $\G\subseteq\Lang$ has \emph{(S)-specification} if there is $t\ge0$ such that for all $w^1,\dots,w^m\in\G$, there are words $v^1,\dots,v^{m-1}$ with $|v^i|=t$ and $w^1v^1w^2\cdots v^{m-1}w^m\in\Lang$. A \emph{decomposition} $\Lang=\Cp\G\Cs$ consists of collections $\Cp,\G,\Cs\subseteq\Lang$ such that every word of $\Lang$ can be written as $uvw$ with $u\in\Cp$, $v\in\G$, $w\in\Cs$. For $M\ge0$, put
\[
\G(M)=\{uvw\in\Lang:\ u\in\Cp,\ v\in\G,\ w\in\Cs,\ |u|\le M,\ |w|\le M\}.
\]
As noted in \cite[\S2.3]{climenhaga2012intrinsic}, the following criterion holds for one-sided as well as two-sided shift spaces; in particular, it applies to the subshifts of $\{0,1\}^{\N}$ considered in this paper. 

\begin{theorem}[{\cite[Thm.~C]{climenhaga2012intrinsic}}]\label{thm:CT}
Suppose that $\Lang$ admits a decomposition $\Lang=\Cp\G\Cs$ such that
\begin{enumerate}[label=\textup{(\Roman*)}]
\item $\G$ has (S)-specification;
\item $h(\Cp\cup\Cs)<h_{\mathrm{top}}(X)$;
\item for every $M$ there is $\tau=\tau(M)$ such that for every $v\in\G(M)$ there are words $u,w$ with $|u|,|w|\le\tau$ and $uvw\in\G$.
\end{enumerate}
Then $X$ is intrinsically ergodic.
\end{theorem}

\subsection{The decomposition}

Fix $k\ge1$ and define
\[
\Cp=\{0^a1^b:\ a,b\ge0,\ b\ge ka\},\qquad
\G=G_k^*,\qquad
\Cs=\{0^a1^b:\ a,b\ge0,\ b\le ka\}.
\]
All three collections are contained in $\Lang_k$: every $0^a1^b$ is a factor of $0^m1^{km}$ for large $m$. Note that a non-empty word of $\Cp$ ends with $1$, a non-empty word of $\Cs$ begins with $0$, and a non-empty word of $\G$ begins with $0$ and ends with $1$.

\begin{remark}\label{rem:natural}
The collection $\Cs$ is exactly the set of prefixes of generators (the empty word included). The set of suffixes of generators is $\{1^b:b\ge0\}\cup\{0^a1^b:\ a\ge1,\ k\mid b,\ b\ge ka\}$, and it is contained in $\Cp$; for $k=1$ the two coincide. For $k\ge2$, the collection $\Cp$ also contains the words $0^a1^b$ with $a\ge1$, $b>ka$ and $k\nmid b$. These are factors lying strictly inside a single generator (for $k=2$, the word $01^3$ inside $0^21^4$); they are neither suffixes nor prefixes of generators, and they cannot be written as $uvw$ with $u$ a suffix of a generator, $v\in G_k^*$ and $w$ a prefix of a generator. This is why the natural decomposition of a coded system into suffixes of generators, concatenations of generators and prefixes of generators \cite{climenhaga2012intrinsic} has to be enlarged here.

The decomposition above is adapted to the systems $X_k$. Note that, in the notation of \cite{climenhaga2012intrinsic}, the collection $\Cs$, whose words come last, consists of \emph{prefixes} of generators, and $\Cp$ of suffixes.  For a general coded system, one would have to add to the suffixes and prefixes of generators the factors lying strictly inside a single generator, and conditions~(II) and~(III) of Theorem~\ref{thm:CT} would have to be checked in each case: these collections may have entropy as large as that of the system itself, and the cost of completing a word of $\G(M)$ to a word of $\G$ may fail to be bounded. Here the collections $\Cp$ and $\Cs$ have at most $n+1$ words of each length $n$, which gives~(II), and the explicit form of $G_k$ gives~(III) (Lemma~\ref{lem:completion}).
\end{remark}

\begin{lemma}\label{lem:htz}
Every $x\in\Lang_k$ containing $10$ can be written uniquely as $x=u\,v\,w$, where $u=0^a1^b$ ends with the $1$ of the first occurrence of $10$ in $x$, $w=0^c1^d$ begins with the $0$ of the last occurrence, and $v\in G_k^*$. Moreover:
\begin{enumerate}[label=\textup{(\roman*)}]
\item $b\ge1$, and if $a\ge1$ then $b=km$ for some $m\ge a$;
\item $c\ge1$, and $d\le kc$; if $d\ge1$, then $0^c1^{kc}$ is the generator containing $w$.
\end{enumerate}
\end{lemma}

\begin{proof}
By Lemma~\ref{lem:boundaries}, the first and last occurrences of $10$ are junctions between generators, so $v$ is a concatenation of whole generators. In (i), if $a\ge1$, then $1^b$ is a run preceded and followed by $0$, hence equal to $1^{km}$ for a generator $0^m1^{km}$ whose zero block contains $0^a$. In (ii), $w$ is a prefix of a generator. If $d\ge1$, its zero block is entirely visible, so the generator is $0^c1^{kc}$ and $d\le kc$.
\end{proof}

\begin{lemma}\label{lem:decomp}
Every word of $\Lang_k$ can be written as $uvw$ with $u\in\Cp$, $v\in\G$ and $w\in\Cs$, and $h(\Cp\cup\Cs)=0$.
\end{lemma}

\begin{proof}
If $x\in\Lang_k$ contains no $10$, then $x=0^a1^b$ lies in $\Cp$ (if $b\ge ka$) or in $\Cs$ (if $b\le ka$). Otherwise, write $x=uvw$ as in Lemma~\ref{lem:htz}: then $u\in\Cp$ by (i), $v\in\G$, and $w\in\Cs$ by (ii). Finally, $\Cp\cup\Cs$ has at most $n+1$ words of length $n$.
\end{proof}

\begin{lemma}[Completion]\label{lem:completion}
Let $M\ge0$. For every $x\in\G(M)$ there are words $p,q$ with $|p|,|q|\le(k+1)M$ such that $pxq\in G_k^*$.
\end{lemma}

\begin{proof}
Write $x=u'v'w'$ with $u'\in\Cp$, $v'\in G_k^*$, $w'\in\Cs$, $|u'|,|w'|\le M$.

\emph{Case 1: $x$ contains no $10$.} If two of $u',v',w'$ were non-empty, a $10$ would appear at their junction, and a $v'$ with at least two generators contains $10$. Hence $x\in G_k\cup\{\eps\}$, in which case we take $p=q=\eps$, or $x\in\{u',w'\}$ with $|x|\le M$. In the latter case, write $x=0^a1^b\neq\eps$, let $m=\max\{a,\lceil b/k\rceil\}\le M$, and take $p=0^{m-a}$, $q=1^{km-b}$. Then $pxq=0^m1^{km}$.

\emph{Case 2: $x$ contains $10$.} Write $x=uvw$ as in Lemma~\ref{lem:htz}, with $u=0^a1^b$ and $w=0^c1^d$. Define
\[
p=\begin{cases}0^{m-a}&\text{if }a\ge1\ (b=km),\\ 0^{r}1^{kr-b}&\text{if }a=0,\ r=\lceil b/k\rceil,\end{cases}
\qquad q=1^{kc-d}.
\]
Then $pu\in G_k$ and $wq\in G_k$, so $pxq=(pu)\,v\,(wq)\in G_k^*$. It remains to bound $p$ and $q$.

If $u'\ne\eps$, then $u'$ ends with $1$ and is followed in $x$ by $0$: otherwise $v'=w'=\eps$ and $x=u'$ contains no $10$. So the first $10$ of $x$ occurs at the junction after $u'$, and $u=u'$. Then $b\le M$, and $|p|\le(k+1)\lceil b/k\rceil\le(k+1)M$. If $u'=\eps$, then $v'\ne\eps$ and the first $10$ of $x$ sits at the end of the first generator of $v'$, so $u\in G_k$ and $p=\eps$.

Symmetrically, if $w'\ne\eps$, then $w=w'$, so $c\le M$ and $|q|\le kM$. If $w'=\eps$, then $w$ is the last generator of $v'$ and $q=\eps$.
\end{proof}

\begin{proposition}[Theorem~\ref{thm:C}(c)]\label{prop:entropy}
$h_{\mathrm{top}}(X_k)=h_{\mathrm{top}}(\Sigma_k)=\frac{\log2}{k+1}$.
\end{proposition}

\begin{proof}
Each generator $0^m1^{km}$ has length $(k+1)m$. Fix $t\ge1$. Each of the $2^{t-1}$ ways of writing $t=m_1+\dots+m_r$ as an ordered sum of positive integers (choose which of the $t-1$ signs in $1+1+\dots+1$ to keep) gives a word
\[
0^{m_1}1^{km_1}\,0^{m_2}1^{km_2}\cdots0^{m_r}1^{km_r}\in G_k^*
\]
of length $(k+1)t$, and different sums give different words, because $m_1,\dots,m_r$ are the lengths of the maximal runs of zeros of the word, read from left to right. Every $w\in G_k^*$ lies in $\Lang(X_k)$, as a prefix of $w^\infty\in X_k$. Hence $\#\Lang_{(k+1)t}(X_k)\ge2^{t-1}$, and since the limit defining $h_{\mathrm{top}}$ exists,
\[
h_{\mathrm{top}}(X_k)=\lim_{t\to\infty}\frac{\log\#\Lang_{(k+1)t}(X_k)}{(k+1)t}\ge\lim_{t\to\infty}\frac{(t-1)\log2}{(k+1)t}=\frac{\log2}{k+1}.
\]
Conversely, every word of $G_k^*$ has the form above, so $\#(G_k^*)_\ell\le2^{\ell/(k+1)}$ for every $\ell\ge0$, and by Lemma~\ref{lem:decomp},
\begin{align*}
\#\Lang_n(X_k) & \le \sum_{i+\ell+j=n}\#\Cp_i\cdot\#(G_k^*)_\ell\cdot\#\Cs_j \\
               & \le \binom{n+2}{2}(n+1)^2\,2^{n/(k+1)} \\
               & \le (n+1)^4\,2^{n/(k+1)}.
\end{align*}
Taking logarithms, dividing by $n$ and letting $n\to\infty$ gives $h_{\mathrm{top}}(X_k)\le\frac{\log2}{k+1}$. The equality $h_{\mathrm{top}}(\Sigma_k)=h_{\mathrm{top}}(X_k)$ is Proposition~\ref{prop:measures}.
\end{proof}

\begin{remark}\label{rem:E}
Let $A_n(p_c)$ be the number of admissible words of length $n$, and let $E(p_c)=\limsup_{n}\frac1n\log A_n(p_c)$ be the entropy of words of \cite[Def.~2]{becker2023symbolic}. Since admissible words belong to $\Lang(\Sigma(p_c))$ (Lemma~\ref{lem:language}), we have $h_{\mathrm{top}}(\Sigma(p_c))\ge E(p_c)$ for every $p_c$. By \cite[Prop.~5]{becker2023symbolic}, $E(\frac{k}{k+1})=\frac{\log2}{k+1}$, so Proposition~\ref{prop:entropy} shows that equality holds for $p_c=\frac{k}{k+1}$. This also sharpens the lower bound $h_{\mathrm{top}}(\Sigma(p_c))\ge(\log2)\,E(p_c)$ of \cite[\S5]{becker2023symbolic}, which is obtained through the Minkowski dimension of an associated subset of $[0,1]$: with the above normalization, that dimension equals $E(p_c)/\log2$, so the inequality it yields is $h_{\mathrm{top}}(\Sigma(p_c))\ge E(p_c)$.
\end{remark}

\begin{theorem}[Theorem~\ref{thm:D}]\label{thm:IE}
For every $k\ge1$, $X_k$ and $\Sigma(\frac{k}{k+1})$ are intrinsically ergodic.
\end{theorem}

\begin{proof}
We verify the hypotheses of Theorem~\ref{thm:CT} for $X_k$ with the decomposition of Lemma~\ref{lem:decomp}.
\begin{enumerate}[label=(\Roman*)]
\item $G_k^*$ is closed under concatenation, so it has (S)-specification with $t=0$. It also has (S)-specification with $t=k+1$, using the connecting word $01^k\in G_k$.
\item $h(\Cp\cup\Cs)=0<\frac{\log2}{k+1}=h_{\mathrm{top}}(X_k)$, by Lemma~\ref{lem:decomp} and Proposition~\ref{prop:entropy}.
\item This is Lemma~\ref{lem:completion}, with $\tau(M)=(k+1)M$.
\end{enumerate}
Hence $X_k$ has a unique MME. By Proposition~\ref{prop:measures}, it is also the unique MME of $\Sigma(\frac{k}{k+1})$.
\end{proof}

Intrinsic ergodicity of $X_k$ can also be deduced from general results. Let $L_{G_k}$ be the subshift of limits of single generators, that is, of sequences all of whose factors are factors of some generator; here $L_{G_k}=\{0^\infty\}\cup\{0^a1^\infty:\ a\ge0\}$, so $h_{\mathrm{top}}(L_{G_k})=0<h_{\mathrm{top}}(X_k)$. By \cite[Thm.~1.4]{pavlov2020coded} (a reformulation of \cite[Thm.~B]{climenhaga2012intrinsic}), $X_k$ has a unique measure of maximal entropy; since $G_k$ is uniquely decipherable (the junctions are the occurrences of $10$), \cite[Thm.~1.5]{pavlov2020coded} and \cite{climenhaga2018towers} give further properties of this measure. We have preferred to verify the hypotheses of Theorem~\ref{thm:CT} directly, because \cite[Thm.~B]{climenhaga2012intrinsic} is formulated in terms of prefixes and suffixes of generators, and the corresponding decomposition in \cite[\S4]{climenhaga2012intrinsic} does not contain the factors lying strictly inside a single generator (Remark~\ref{rem:natural}).

%==================================================================
\section{Questions}\label{sec:questions}
%==================================================================

\begin{question}
For which $p_c\in(0,1)$ is $\Sigma(p_c)$ transitive? So far we know the answer only for $p_c=\frac{k}{k+1}$: yes for $k=1$, and no for $k\ge2$ (Theorem~\ref{thm:C}).
\end{question}

\begin{question}
Is $\Sigma(p_c)$ intrinsically ergodic for every $p_c\in(0,1)$? 
\end{question}

%%=====================
\section*{Acknowledgments}

The author is grateful to Professor Rodrigo Bissacot for many helpful discussions and valuable comments, and for his support during the preparation and editing of this manuscript. This study was financed, in part, by the São Paulo Research Foundation (FAPESP), Brasil, Process Number \#2025/07015-0. The author used Claude Opus 5.5 (Anthropic) to improve the English language of this text.
%==================================================================
\bibliographystyle{amsplain}
\bibliography{bibliography}

\end{document}